\documentclass[reqno]{amsart}
\usepackage{graphicx}
\usepackage{amsmath,amssymb}

\newtheorem{thm}{Theorem}[section]

\newtheorem{lem}[thm]{Lemma}

\begin{document}
\large 
  \begin{center}
  \textbf{An elementary construction of the Wiener measure}\\[0.6 cm]
  R. P. Pakshirajan* and M. Sreehari**\\[0.5 cm] \end{center}
  \begin{footnotesize}
\makeatletter{\renewcommand*{\@makefnmark}{}
\footnotetext{*R. P. Pakshirajan, 227, 18th Main, 6th Block, Koramangala, Bengaluru -560095, India.\\
Corresponding author:** M. Sreehari, 6-B, Vrundavan Park, New Sama Road, Vadodara, 390024, India.\\
E-\textit{mail addresses}: vainatheyarajan@yahoo.in (R. P. Pakshirajan);  msreehari03@yahoo.co.uk (M. Sreehari)  }}
\end {footnotesize}
 \textbf{Abstract.}\\
 Our construction of the Wiener measure on $\mathfrak{C}$ consists in first defining a  set function $\varphi$\ on the class of all compact sets based on certain $n$-dimensional normal distributions, $n = 1,\ 2,\ldots$\ using the structural relation at (1) below. This structural relation, discovered by the first author, is recorded in his book [2] on page 130. We then define a measure $\mu$ on the  Borel $\sigma$-field of subsets of $\textbf{C}$ which is the Wiener measure. 
 \section {Introduction and preliminary results}
 Let $\textbf{C} = \textbf{C}[0,\ 1]$\ denote the space of real valued continuous functions defined on $[0,\ 1]$, all vanishing at $0$\ and endowed with
 the norm $\Arrowvert x \Arrowvert = \sup\limits_{0 \le t \le 1}|x(t)|,\ x \in \textbf{C}$. Write $\rho(x,\ y) = \Arrowvert x - y \Arrowvert,\ x,\ y \in \textbf{C}$. Let $\mathfrak{C}$\ denote the Borel
 $\sigma$-field (i.e., the $\sigma$-field generated by the open subsets) of \textbf{C}. For $t_j \in (0,\
 1],\ 1 \le j \le k$, define projection operators $\pi_{t_1,\ t_2,\ldots,\ t_k}$\ thus. For $x \in
 \textbf{C},\ \pi_{t_1,\ t_2,\ldots,\ t_k}(x) = \big(x(t_1),\ x(t_2),\ldots,\ x(t_k)\big)$. This maps
 \textbf{C} into $R^k$. Assume $R^k$\ is endowed with the usual metric and denote the resulting Borel
 $\sigma$-field by $\mathfrak{R}^k$. We note $\pi_{t_1,\ t_2,\dots.,\ t_k}$\ is a continuous map and
 $\mathfrak{C} $\ measurable.\\\\ Let $\nu_{t_1,\ t_2,\ldots,\ t_k}$\ denote the measure
 generated on $\mathfrak{R}^k$\ by a $k$-variate normal distribution whose components have zero means and
 the covariance between the $\text{j}^{\text{th}}$\ and the $\text{r}^{\text{th}}$\ components is
 $\min\{t_j,\ t_r\}$. We note the variance covariance matrix is positive definite and hence the
 distribution is non-singular. Denote the corresponding frequency function by $f_{t_1,\ t_2,\ldots.,\ t_k}$. Let $\alpha_{t_1,\ t_2,\ldots,\ t_k}$\ denote the measure generated on the sub $\sigma$-field
 $\pi^{-1}_{t_1,\ t_2,\ldots,\ t_k}(\mathfrak{R}^k)$\ by the mapping $\pi_{t_1,\ t_2,\ldots,\ t_k}$. Note $\alpha_{t_1, t_2, \ldots, t_k} (A)$ decreases as $k$ increases for $ A \in \pi_{t_1, t_2,\ldots, t_k}^{-1} (\mathcal {R}^k)$.\\\\
  For the usual method of construction of the Brownian motion process / Wiener measure, refer to chapter 2 in  Karatzas and Shreve [1].\\\\ 
 Let $T =\{t_1, t_2,\ldots\}$\ be a countable dense subset of $[0,\ 1]$. Let $K$\ be a compact
 subset of $\textbf{C}$. Hence $\pi_{t_1, t_2,\ldots, t_n}(K)$\ is a compact subset of $R^n$. Then the following structural relation holds:
 (ref. pp 130-131 in Pakshirajan [2])
 \begin{equation} \label{E1}
 K = \bigcap_{n = 1}^\infty \pi^{-1}_{t_1,t_2,\ldots,t_n}\pi_{t_1, t_2,\ldots,t_n}(K).	
 \end{equation}
For completeness we include  proof of (\ref{E1}).\\
Since $ K \subset \pi^{-1}_{t_1, t_2,\ldots,t_n}\pi_{t_1, t_2,\ldots,t_n}(K) $, it follows that $K$ is a subset of the set on the right side of (\ref{E1}).\\
Let now $x$ be an arbitrary member of the set on the right side of (\ref{E1}). Hence for every $n $, $x\in \pi^{-1}_{t_1, t_2,\ldots,t_n}\pi_{t_1, t_2,\ldots,t_n}(K)$. Then there exists $y_n \in K$ such that $\pi_{t_1,\ t_2,\ ....,\ t_n} x= \pi_{t_1,\ t_2,\ ....,\ t_n} y_n$. Since $K$ is a compact set the sequence $(y_n)$ admits a convergent subsequence, say, ($y_m$) converging to, say, $y_0 \in K$ in the metric $\rho$. This implies $y_m(t) \rightarrow y_0(t)$ for all $t \in [0, 1]$. Fix $r$. Then for all $m > r$, $x(t_r)=y_m(t_r)$. Letting $m \rightarrow \infty$ we observe $x(t_r) =y_0(t_r)$. Since this is true for every $t_r \in T$, since $T$ is a dense subset of $[0, 1]$ and since  the functions are continuous, we get $x\equiv y_0$. Thus $x \in K$ and the result follows.

\begin{thm}
 $\alpha_{t_1,\ t_2,\ldots.,\ t_n}(\pi^{-1}_{t_1,\ t_2,\ldots,\ t_n}\pi_{t_1,\ t_2,\ldots,\ t_n}(K)),\ n = 1,\ 2,\ldots.$\ is a
monotonic decreasing sequence of numbers. 

\end{thm}
\noindent \textbf{ Proof.}
\begin{equation*}
	\alpha_{t_1, t_2,\ldots, t_{n +1}}\big(\pi^{-1}_{t_1,t_2,\ldots, t_{n + 1}}\pi_{t_1,t_2,\ldots, t_{n + 1}}(K)\big)
 \end{equation*}
\begin{equation*}
 =\int_{\pi_{t_1, t_2,\ldots, t_{n + 1}}(K)} f_{t_1, t_2,\ldots, t_{n + 1}}(u_1, u_2,\ldots, u_{n + 1}) du_1 d u_2\dots,du_{n + 1}
  \end{equation*}
\begin{equation*}
 \le \int_{\pi_{t_1, t_2,\ldots,t_{n}(K) \times (- \infty, \infty)}} f_{t_1, t_2,\ldots, t_{n + 1}}(u_1,u_2,\ldots,\ u_{n + 1})du_1 du_2\ldots,du_{n + 1}
 \end{equation*}
\begin{equation*}
 \le \int_{\pi_{t_1,t_2,\ldots, t_n}(K)} f_{t_1,t_2,\ldots, t_n}(u_1, u_2,\ldots, u_n)du_1 du_2\ldots du_{n} 
 \end{equation*}
\begin{equation*}
 \le\alpha_{t_1, t_2,\ldots, t_{n}}\big(\pi^{-1}_{t_1, t_2,\ldots, t_{n}}\pi_{t_1, t_2,\ldots, t_{n}}(K))
\end{equation*}
completing the proof.\\
Define set function $\varphi$\ on the compact sets $K$\ of \textbf{C} :
\begin{equation}\label{E2}
	\varphi(K) = \lim_{n \rightarrow \infty}\alpha_{t_1,\ t_2,\ ....,\ t_{n}}\big(\pi^{-1}_{t_1,\ t_2,\ ....,\ t_{n}}\pi_{t_1,\ t_2,\ ....,\ t_{n }}(K)\big).
	\end{equation}
Since $\pi_{t_1,\ t_2,\ ....,\ t_{n}}(\textbf{C}) = R^n$\ and since (1) is true with $K$\ replaced by $\textbf{C}$,  we set
$\varphi(\textbf{C}) =1$. 

We shall now discuss some properties of the set function $\varphi$.
\begin{lem}\label{L2}
Let $K_1,\ K_2$\ be two disjoint compact subsets of $\textbf{C}$. Then\\ $\varphi(K_1 \cup K_2) = \varphi(K_1) + \varphi(K_2).$
\end{lem}
\noindent \textbf{Proof.}\\
Since $K_1 \cup K_2$\ is a compact set
\begin{equation*}
	\varphi(K_1 \cup K_2) = \lim_{n \rightarrow \infty}\alpha_{t_1,t_2,\ldots,		t_{n}}\big(\pi^{-1}_{t_1, t_2,\ldots,  t_{n}}\pi_{t_1, t_2,\ldots,  t_{n }}(K_1 \cup
	K_2)\big)
\end{equation*}
\begin{equation*}
=	\lim_{n \rightarrow \infty}\alpha_{t_1, t_2,\dots,
		t_{n}}\big(\pi^{-1}_{t_1, t_2,\ldots,  t_{n}}\{\pi_{t_1, t_2,\dots,  t_{n }}(K_1) \cup\\
	\pi_{t_1, t_2,\dots, t_{n }}(K_2)\}\big)
\end{equation*}
\begin{equation*}
	= \lim_{n \rightarrow
		\infty}\Big[\alpha_{t_1, t_2,\dots,  t_{n}}(\pi^{-1}_{t_1, t_2,\ldots,  t_{n}}\{\pi_{t_1, t_2,\ldots,  t_{n }}(K_1)\}) +
\end{equation*}
\begin{equation*}	
	+ \alpha_{t_1, t_2,\ldots,  t_{n}}(\pi^{-1}_{t_1, t_2,\ldots, 
		t_{n}}\{\pi_{t_1, t_2,\ldots,  t_{n }}(K_2)\})\Big]
\end{equation*}
since the events are mutually exclusive and since $\alpha_{t_1,t_2,\ldots,  t_{n}}$ is a measure.  Then 
\begin{equation*}
\varphi(K_1 \cup K_2) =	\lim_{n \rightarrow \infty}\alpha_{t_1, t_2,\ldots,t_{n}}(\pi^{-1}_{t_1, t_2,\ldots,t_{n}}\{\pi_{t_1, t_2,\ldots, t_{n }}(K_1)\})+
\end{equation*}
\begin{equation*}
+ \lim_{n \rightarrow\infty}
	\alpha_{t_1, t_2,\ldots,t_{n}}(\pi^{-1}_{t_1, t_2,\ldots,  t_{n}}\{\pi_{t_1, t_2,\ldots, t_{n }}(K_2)\})
\end{equation*}
since the individual limits exist. We thus have  $\varphi(K_1 \cup K_2) = \varphi(K_1) + \varphi(K_2).$\\\\
\textbf{Remark 1.} We have the following observations from the earlier discussion:\\
a) $\varphi$ is finitely additive on the collection of compact sets.\\
b) $\varphi (K) \geq 0$ for all compact sets $K$.\\
c) If $K_1, K_2$ are compact sets and $K_1 \subset K_2$ then from (\ref{E2}),  $\varphi(K_1) \leq \varphi(K_2)$.\\
We now discuss some limiting properties of $\varphi(K_n)$.\\\\
\noindent \textbf{Definition}. We call a set, in a topological space, a boundary set if it is a closed set with a null interior. The boundary of  a set $A$ (i.e.,   $\overline{A} \sim\textit{Int}\;A$ ) will be denoted by  $\partial A$. \\
We note that the boundary of a set is a boundary set.
\begin{lem} \label {L4}
 Let $f$\ be a continuous map of a topological space $X$\ into a topological space $Y$. If $A \subset X$\ is a boundary set and if $f(A)$ is a closed set, then  $f(A)$ is a boundary set \ in $Y$. 
 \end{lem}
\noindent \textbf{Proof.}\\
If $f(A)$\ contains an open set $G$, then it is clear that $f^{-1}(G)$\ is an open subset of $A$, which is not possible since it is a boundary set.
\begin{thm} \label{T2}
	(i) If $E \subset \textbf{C}$\ is a compact boundary set, then $\varphi(E) = 0$.\\
	(ii) Let $K_n,\ n \ge 1$\ be compact boundary sets and let $K_n \rightarrow K$\ where $K$\ is a compact set. Then $\varphi(K) = 0$.\\
	(iii) Let $K_n, n \geq 1$, be compact sets such that $K_n \downarrow K$. Then $\varphi(K_n) \downarrow \varphi(K).$\\
	\end{thm}
\noindent \textbf{Proof.}\\
(i) The compact set $\pi_{t_1, t_2,\ldots, t_r}E$\ is, by  Lemma \ref{L4}, a boundary set in $R^r$. Hence its  Lebesgue measure is zero. Since the measure generated by the multivariate normal distribution is absolutely continuous with respect to the Lebesgue measure, it follows that $\alpha_{t_1,t_2,\ldots, t_r}(\pi^{-1}_{t_1,t_2,\ldots, t_r}\pi_{t_1,t_2,\ldots, t_r}E) = 0$. Since this is true for every $r \ge 1$, we conclude $\varphi(E) = 0$. \\
(ii) We note that, by part (i), $\varphi(K_n) = 0,\ n \ge 1$. Further for all $r \ge 1,\ n \ge 1,\ \alpha_{t_1,t_2,\ldots, t_r}(\pi^{-1}_{t_1,t_2,\ldots, t_r}\pi_{t_1,t_2,\ldots, t_r}K_n) = 0$\ by an appeal to the Lemma \ref{L4} and arguing as in Part (i).  If possible let $\varphi(K) > 0$. Hence for all $r \ge 1$,\\ $\alpha_{t_1,t_2,\ldots,  t_r}(\pi^{-1}_{t_1,t_2,\ldots, t_r}\pi_{t_1,t_2,\ldots, t_r}K) > 0$. Since the convergence \\ $\alpha_{t_1,\ldots,  t_r}(\pi^{-1}_{t_1,t_2,\ldots, t_r}\pi_{t_1,t_2,\ldots, t_r}K_n) \rightarrow \alpha_{t_1,t_2,\ldots,  t_r}(\pi^{-1}_{t_1,t_2,\ldots, t_r}\pi_{t_1,t_2,\ldots, t_r}K) $\ holds and since this leads to the absurd result that a sequence of zeros converges to a positive number, it follows that $\varphi(K) = 0$.\\
(iii) Since the sequence $(\varphi(K_n))$\ is monotonic decreasing, it is enough to show that, given $\varepsilon > 0$, there exists $K_N$\ such that $\varphi(K_N) < \varphi(K) + \varepsilon$.\\
We note $K$\ is compact. Hence given $\varepsilon > 0$, we can find $r \ge 1$\ such that
\begin{equation} \label{E3}
	\varphi(K) > \alpha_{t_1, t_2,\ldots, t_{\ell}}\big(\pi^{-1}_{t_1, t_2,\ldots, 
		t_{\ell}}\pi_{t_1, t_2,\ldots, t_{\ell}}(K)\big) - \varepsilon\
\end{equation}
for all $\ell \ge r.$ Since $K_n \downarrow K$, for all $\ell \geq 1$ we have
\begin{equation*}
	\pi^{-1}_{t_1, t_2,\ldots, t_{\ell}}\pi_{t_1, t_2,\ldots, t_{\ell}}(K_n) \downarrow \pi^{-1}_{t_1, t_2,\ldots, t_{\ell}}\pi_{t_1, t_2,\ldots, t_{\ell}}(K).
\end{equation*}
For fixed $\ell$ we then have, as $n\rightarrow \infty$
\begin{equation*}
	\alpha_{t_1, t_2,\ldots, t_{\ell}}\big(\pi^{-1}_{t_1, t_2,\ldots,
		t_{\ell}}\pi_{t_1, t_2,\ldots, t_{\ell}}(K_n)\big) \downarrow \alpha_{t_1, t_2,\ldots, t_{\ell}}\big(\pi^{-1}_{t_1, t_2,\ldots, t_{\ell}}\pi_{t_1, t_2,\ldots, t_{\ell}}(K)\big).
\end{equation*}
Take $\ell = r$. We can find $N = N(r)$\ large such that
\begin{equation*} 
	\alpha_{t_1, t_2,\ldots, t_{r}}\big(\pi^{-1}_{t_1, t_2,\ldots,  t_{r}} \pi_{t_1, t_2,\ldots, t_{r}}(K_N)\big) < \alpha_{t_1, t_2,\ldots, t_{r}}\big(\pi^{-1}_{t_1,\ t_2,\ldots, t_{r}}\pi_{t_1, t_2,\ldots, t_{r}}(K)\big) + \varepsilon.
\end{equation*}
This, together with (\ref{E3}), yields
\begin{equation*}
	\varphi(K) + \varepsilon > \alpha_{t_1, t_2,\ldots, t_{r}}\big(\pi^{-1}_{t_1, t_2,\ldots, t_{r}} \pi_{t_1, t_2,\ldots, t_{r}}(K_N)\big) - \varepsilon > \varphi(K_N) - \varepsilon.
\end{equation*}
Since $\varepsilon > 0$\ is arbitrary, the claim follows.\\
\section{The Wiener measure}
In this Section we introduce a new set function in terms of $\varphi$ on the Borel $\sigma$-field $\mathfrak{C}$ and study its properties to show that it is indeed the Wiener measure.\\
For  arbitrary measurable sets  $A \in \mathfrak{C}$\, define
\begin{equation} \label {E4}
	\mu(A) = \sup\limits_{K \subset A,\ K\ \text{compact}}\varphi(K).
\end{equation}
At the outset we observe that for compact sets $K$, $\mu(K)=\varphi(K)$ and hence all the properties noted in the previous Section for $\varphi$ also hold for $\mu$. Further the definition 
implies (i) that if $A \subset B, A, B \in \mathfrak{C}\ \text{then}\ \mu(A) \le \mu(B)$\ and (ii) that there exists an increasing sequence ($K_n$) of compact sets, $K_n \subset A$\ such that $\mu(A) = \lim\limits_{n \rightarrow \infty}\mu(K_n)$.\\
The sets $K_n$ can be chosen to be monotonic increasing.\\
\textbf{Remark 2.}  This does not mean that $K_n \uparrow A$. i.e., $\bigcup\limits_{n = 1}^{\infty}K_n$\ can be a proper subset of $A$. 
To see this, take $v \in \textbf{C}, \|v\|=1.$ Let $K_n=\{\lambda v,\; 0 \leq \lambda \leq 1-\frac{1}{n}\}$ and  $A=\{\lambda v, \; 0\leq \lambda \leq 1 \}.$  However, if  $K_n=\{\lambda v,\; 0 \leq \lambda \leq 1-\frac{1}{n}\}\cup \{v\}$ then both $K_n$ and $A$ are compact and  $K_n\uparrow A$.\\\\
We next discuss further properties of $\mu$ that enable us to claim that $\mu$ is indeed a probability measure.
\begin{lem} \label{L3}
	Let $ A,\ B \in \mathfrak {C},\ A \cap B = \emptyset$. Then $\mu(A \cup B) = \mu(A) + \mu(B)$.
\end{lem}
\noindent \textbf {Proof}.\\ Let $E \subset A,\ F \subset B$\ be compact sets. We have, from Lemma \ref{L2}\\
$\mu(E) + \mu(F) = \mu(E \cup F) \le \mu(A \cup B)$.\\
Thus $\mu(A \cup B) \ge \mu(A) + \mu(B)$. It remains to be shown that $\mu(A \cup B) \le \mu(A) + \mu(B).$\\
Given $\varepsilon > 0$, we can find a compact set $K,\;\ K \subset A \cup B$\ such that\\ 
$\mu(A \cup\ B) - \varepsilon < \mu(K)$.\\
Case 1. The distance  $d(A,\ B) = q > 0$.\\ Consider an arbitrary sequence $(x_n)\ \text{in}\ K \cap A$. Since it is a sequence in $K$, it contains a convergent subsequence, converging to, say, $x_0$. This $x_0$\ has to be in $K \cap A$\ or in $K \cap B$. Since the sequence lies in $K \cap A$\ and since $d(K \cap A,\ K \cap B) \ge q > 0$, we conclude $x_0 \in K \cap A$. Thus we see that every sequence in $K \cap A$\   
contains a convergent subsequence converging to a point in $K \cap A$. This means $K \cap A$\ is a compact set. Similarly, $K \cap B$\ is a compact set. Summarising, we conclude that every compact subset of $A \cup B$\ is the union of a compact subset $E$\ of $A$\ and a compact subset $F$\ of $B$. We get $\mu(A \cup\ B) - \varepsilon < \mu(K) = \mu(E \cup F) = \mu(E) + \mu(F) \le \mu(A) + \mu(B)$.  That $\mu(A \cup B) \le \mu(A) + \mu(B)$\ is now immediate.\\ 
Case 2. $d(A,\ B) = 0$.\\ This case assumption implies that $Q = \bar{A} \cap \bar{B} \neq \emptyset$. Again in this case one or both the sets $K \cap A,\ K \cap B$\ can fail to be compact. Since the other case admits to being similarly argued, let us assume that neither of the two sets is compact. 
$K \subset A \cup B$\ can not be compact if any convergent sequence in it converges to a point outside $K$. i.e., if convergent sequences in $K \cap A$\ or in $K \cap B$\ converge to points outside these sets. Thus $K$\  can be a compact subset only if $E = K \cap A$\ and $F = K \cap B$\ are compact. And the arguments and the conclusion in case 1 hold.\\ With this the proof is complete.\\
\textbf{Remark 3.}\\Immediate consequences of Lemma \ref {L3} are :\\
a) If $A_k,\ 1 \le k \le n$\ is a collection of $n$\ mutually exclusive sets, then $\mu(\bigcup\limits_ {k = 1}^nA_k) = \sum\limits_{k = 1}^n\mu(A_k)$\ and\\
b) If $A,\ B \in \mathfrak{C}$\ and if $A \subset B$, then $\mu(B \sim A) = \mu(B) - \mu(A)$.
\noindent \begin {thm} \label{T3}
	(i) Let $K_n \uparrow,\; K_n$ be a compact boundary set for $ n\geq 1$. If $A=\bigcup_{n=1}^\infty K_n$, then $\mu(A) = 0.$\\
	 (ii) Suppose $K, K_n, n \ge 1$ are compact subsets with $K_n \uparrow K$.  Then $\mu(K_n) \uparrow \mu(K)$.\\
	(iii) Let $K_n \uparrow,\ K_n$\ be compact. Let $A = \bigcup_{n = 1}^{\infty}K_n$. Then $\mu(A) = \lim_{n \rightarrow \infty}\mu(K_n)$.\\
\end{thm}
\noindent \textbf{Proof.}\\
(i) If $\mu (A) > 0$ then given $\varepsilon > 0$ we can find a compact set $K \subset A$ such that $\mu (A) < \mu (K) +\varepsilon.$ We note $\mu(K\bigcap K_n)=0$ and $K\bigcap K_n \uparrow K$. Now appeal to Theorem \ref{T2} (ii) and conclude $\mu(K) =0.$ This implies $\mu(A) < \varepsilon.$ Since $\varepsilon$ is arbitrary it follows $\mu(A) =0.$ \\
\noindent(ii) We are to prove $\mu(K) - \mu(K_n) \rightarrow 0$. Write $A_n = K \sim K_n \downarrow \emptyset$. Interior $B_n$\ of $A_n$\ is $(Int K) \cap (K \sim K_n)$, since $K_n^{\prime}$\ is an open set. We note that $\bar{A}_n = \bar{B}_n = A_n \cup C_n$\ where $	C_n = \partial K_n$.\\
We note that $A_n \cap C_n = \emptyset$, that $\bar{A}_n$\ is a compact set, that the ($C_n$) is a sequence of disjoint compact boundary sets and that $(\bar{A}_n)$\ is a monotonic decreasing sequence, the last assertion being true since $A_n \downarrow$. We will show $\mu(\bar{A}_n) \rightarrow 0$\ and that will complete the proof of this part.\\ 
Now, if $A = \bigcap\limits_{n = 1}^{\infty}\bar{A}_n$, then by Theorem \ref{T2}(iii), $\lim_{n \rightarrow\infty}\mu (\bar{A}_n) = \mu(A)$. \\ Let $D = \bigcup\limits_{n = 1}^{\infty}C_n$.
We note that $\bar{A}_n = A_n \cup C_n \subset A_n \cup D$. Hence $A \subset (\bigcap\limits_{n = 1}^{\infty}A_n) \cup D = D$, leading to $\mu(A) \le \mu(D)$.\\ We note that for each $n,\ \mu(C_n) = 0$\ since $C_n$\ is a compact boundary set. (ref. Theorem \ref{T2}, part (i) )
If $D_n = \bigcup_{j = 1}^nC_j$\ then, by part (a) of Remark 1, $\mu(D_n) = \sum\limits_{j = 1}^n \mu(C_j) = 0$. $D_n$, being the finite union of disjoint compact boundary sets, is a compact boundary set. Since $D_n \uparrow D$, it follows by part(i) above 
that  $\mu(D) = 0$. Thus $\mu(A) = 0$. With this the proof of (ii) is complete.\\
\noindent(iii) The hypothesis implies 
\begin{equation} \label{E5}
	\lim_{n \rightarrow \infty}\mu(K_n) \le \mu(A).
\end{equation}
If $A$\ is a compact set, then the claim follows from part (ii). If $A$\ is not a compact set, then, given $\varepsilon$, a compact set $K \subset A$\ can be found such that $\mu(K) > \mu(A) - \varepsilon$. Now, by part (ii), $\lim\limits_{n \rightarrow \infty}\mu(K \cap K_n) = \mu(K) > \mu(A) - \varepsilon$. Hence $\lim_{n \rightarrow \infty}\mu(K_n) \ge \lim_{n \rightarrow \infty}\mu(K \cap K_n) \ge \mu(A) - \varepsilon$\\ This, together with (5) and use of the fact that $\varepsilon >0$\ is arbitrary completes the proof of (iii).\\
\textbf{Remark 4}\\ The result (ii) above is for compact sets but we could not get it in Section 1 because its proof needs the result in (i).\\
We now extend some of the above results (proved for compact sets) to measurable sets.
\begin{thm} \label {T4}
	(i) Let $ A_n \in \mathfrak{C},\ n \ge 1;\ K$\ be compact. Let $A_n \uparrow K$. Then $\mu(A_n) \uparrow \mu(K)$.\\
(ii) Let $A_n,\ A \in \mathfrak{C}$, with $A_n \uparrow A$. Then $\mu(A_n) \uparrow \mu(A)$.
\end{thm}
\noindent\textbf{Proof.}\\

\noindent(i) The hypothesis $(A_n \uparrow K)\ \Longleftrightarrow (E_n = K \sim A_n \downarrow \emptyset )$. Note that 
$\mu(A_n) \uparrow.$  We note $\bar{E}_n \subset K$\ is a compact set for each $n$\ and $\bar{E}_n \downarrow$. Now,
\begin{equation*}
	\bigcap_{n = 1}^{\infty}\bar{E}_n = \bigcap_{n = 1}^{\infty}\big(E_n \cup \partial E_n\big) \subset \bigcap_{n = 1}^{\infty}(E_n \cup D) = D
\end{equation*}
where $ D = \bigcup_{n = 1}^{\infty}\partial E_n$  since $\bigcap\limits_{n = 1}^{\infty}E_n = \emptyset$. Theorem \ref{T2} (iii) applies and we get  $\lim_{n \rightarrow \infty}\mu(\bar{E}_n) = \mu\big(\bigcap_{n = 1}^{\infty}\bar{E}_n\big) \le \mu(D)$. Since $D$\ is the union of a countably many disjoint  compact boundary sets, it follows from Theorem \ref{T2}(ii) that $\mu(D) = 0$. Thus $\lim_{n \rightarrow \infty}\mu(\bar{E}_n) = 0$. Since $\partial E_n$\ is a boundary set, its $\mu$-measure is zero. Hence $\mu(E_n) \rightarrow 0$. i.e., $\mu(K) - \mu(A_n) \rightarrow 0$, as was to be shown.\\
\noindent(ii)  Given $\varepsilon > 0$\ find compact set $K \subset A$\ such that $\mu(A) - \varepsilon < \mu(K)$. Since $K \cap A_n \uparrow K$, part (i) applies and we get
\begin{equation*}
	\mu(K) = \lim_{n \rightarrow \infty}\mu(K \cap A_n) \le \lim_{n \rightarrow \infty}\mu(A_n).
\end{equation*}
Thus $\lim_{n \rightarrow \infty}\mu(A_n) \ge \mu(K) \ge \mu(A) - \varepsilon$, leading to $\lim_{n \rightarrow \infty}\mu(A_n) \ge \mu(A)$. This together with the obvious inequality that $\lim_{n \rightarrow \infty}\mu(A_n) \le \mu(A)$  completes the proof.
\begin{thm} \label{T5}
	$\mu$\ defined as above at (\ref{E4}) on $\mathfrak{C}$\ is a probability measure.
\end{thm}	
\noindent\textbf{Proof.}\\
Let $A_n \in \mathfrak{C},\ n \ge 1$\ be a sequence of mutually exclusive events. Let $A = \bigcup_{n = 1}^{\infty}A_n = \bigcup_{n = 1}^{\infty}B_n$\ where $B_n = \bigcup_{k = 1}^nA_k$. Since $B_n \uparrow A$, Theorem  \ref{T4} (i) applies and then we have by the observations following (\ref{E4})
\begin{equation*}
	\mu(A) = \lim_{n \rightarrow \infty}\mu(B_n) = \lim_{n \rightarrow \infty}\sum_{k = 1}^n\mu(A_k)= \sum\limits_{k = 1}^{\infty}\mu(A_k).
\end{equation*}
i.e., $\mu$\ on $\mathfrak{C}$\ is countably additive. Since $\mu (A) \geq 0$ for $A\in \mathfrak{C} $ and $\mu(\textbf{C})=1$ it follows that $\mu$ is a probability measure.\\
\noindent \textbf{Remark 5.}\\ (i) Let $\{u_1, u_2, \ldots, u_p\} \subset [0, 1]$ be arbitrary. We can find $\{t_{1,n}, t_{2, n}, \ldots, t_{p,n}\} \subset T$ such that $t_{j, n}\rightarrow u_j, \; 1 \leq j \leq p$. Then we will have  $\pi_{t_{1, n}, t_{2, n}, \ldots, t_{p, n}}\; x\rightarrow \pi_{u_1, u_2, \ldots, u_p}\;x$ for every $x \in \bf{C}$. Hence the random vector  $\pi_{t_{1, n}, t_{2, n}, \ldots, t_{p, n}}\;\stackrel{D}{\rightarrow}  \pi_{u_1, u_2, \ldots, u_p}\;$. Using characteristic functions, we note $\pi_{u_1, u_2, \ldots, u_p}$ has, under $\mu$,  a multi-
variate normal distribution with $cov (\pi_{u_j}, \pi_{u_r})= \min\{u_j, u_r\}$. It follows
now that $\mu$ is the Wiener probability measure.\\
(ii) The co-ordinate process $\{\pi_t,\ t \ge 0\}$\ is known as the Brownian motion process. \\\\
\noindent\textbf{References}\\\\
1. Karatzas, I. Shreve, S. E.: Brownian motion and stochastic calculus, Springer  (1988).\\
2. Pakshirajan, R. P.:  Probability Theory (A foundational course), Hindustan Book Agency, New Delhi (2013).

\end{document}